\documentclass[11pt]{article}
\usepackage{enumerate}
\usepackage{amsmath,amsthm}
\usepackage{amssymb, latexsym}
\usepackage[mathscr]{euscript}
\usepackage{color}
\usepackage{array}
\usepackage{flafter}
\usepackage{layout}
\usepackage{graphicx}
\usepackage[applemac]{inputenc}
\usepackage{hyperref}
\theoremstyle{definition}
  \newtheorem{defi}{Definition}[section]
  
\theoremstyle{plain}
  \newtheorem{teo}{Theorem}[section]
  
  \newtheorem{lema}{Lemma}[section]
  \newtheorem{prop}{Proposition}[section]

\title{Percolation for the stable marriage of Poisson and Lebesgue with random appetites}
\author{Daniel Andr\'es D\'{\i}az-Pach\'on\thanks{Biostatistics Division - University of Miami. DDiaz3@miami.edu}} 
\date{\today}
\allowdisplaybreaks
\begin{document}
\maketitle
\begin{abstract}
\noindent Let $\Xi$ be a set of centers chosen according to a Poisson point process in $\mathbb R^d$. Consider the allocation of $\mathbb R^d$ to $\Xi$ which is stable in the sense of the Gale-Shapley marriage problem, with the additional feature that every center $\xi\in\Xi$ has a random appetite $\alpha V$, where $\alpha$ is a nonnegative scale constant and $V$ is a nonnegative random variable. Generalizing previous results by Freire, Popov and Vachkovskaia (\cite{FPV}), we show the absence of percolation when $\alpha$ is small enough, depending on certain characteristics of the moment of $V$.
\\
\textbf{Key words}:  Poisson process, percolation, stable marriage, random appetite.
\end{abstract}

\section{Introduction}

Let $\Xi$ be a set of \textit{centers} chosen according to an homogeneous Poisson point process of law $\textbf P$ with intensity $\lambda>0$. Let $V$ be a nonnegative random variable (r.v.)\ with law $\mathbb P$. We write  $\textbf E$ and $\mathbb E$ for the expectations under $\textbf P$ and $\mathbb P$, respectively. Consider also a constant $\alpha\geq0$. 

The product $\alpha V$ will be called (random) \textit{appetite}. Every center $\xi\in\Xi$ will have an appetite in such a way that appetites corresponding to different centers are independent of each other and identically distributed. The appetites are also independent of $\Xi$.  Keeping in mind the independence (with a little abuse of notation), we will denote the appetite of a center $\xi$ by $\alpha V_\xi$, when it does not lead to confusion. 

The elements $x$ of $\mathbb R^d$ will be called \textit{sites}. Let $\psi$ be a function of $\mathbb R^d$ to $\Xi\cup\{\infty,\Delta\}$; we will call $\psi$ an \textit{allocation}. Let $\mathcal L$ be the Lebesgue measure. Given a center $\xi$, its \textit{territory} will be given by $\psi^{-1}(\xi)$. The allocation $\psi$ satisfies the restriction $\mathcal L\left[\psi^{-1}(\xi)\right]\leq\alpha V_\xi$, for all $\xi\in\Xi$ (i.e., the volume of a center's territory is bounded by its appetite). A center $\xi$ will be called \textit{unsated} if $\mathcal L\left[\psi^{-1}(\xi)\right]<\alpha V_\xi$; we will say that $\xi$ is \textit{sated} if $\mathcal L\left[\psi^{-1}(\xi)\right]=\alpha V_\xi$. If $\psi(x)=\xi$, for some $\xi\in\Xi$, we will say that $x$ is \textit{claimed}; if $\psi(x)=\infty$, we will call it \textit{unclaimed}. Finally, the set $\Delta$ is a convenience to denote a $\mathcal L$-null set of sites not allocated to $\Xi\cup\{\infty\}$; that is, $\mathcal L[\psi^{-1}(\Delta)]=0$.

Consider a site $x$ and a center $\xi$ such that $\psi(x)\notin\{\xi,\Delta\}$. We will say that the site $x$ \textit{desires} the center $\xi$ if $|x-\xi|<|x-\psi(x)|$ or if $\psi(x)=\infty$ (here and henceforth $|\cdot|$ is the Euclidean norm). If there exists $x'\in\psi^{-1}(\xi)$ such that $|x-\psi(\xi)|<|x'-\psi(\xi)|$, or if $\xi$ is unsated, we will say that the center $\xi$ \textit{covets} the site $x$. A pair $(x,\xi)$ will be \textit{unstable} if $x$ desires $\xi$ and $\xi$ covets $x$. If $\psi$ does not produce any unstable pair it will be called \textit{stable}. Finally, call $\mathcal C$ the closure of the set of claimed sites.

Hoffman, Holroyd, and Peres gave an explicit construction of the stable function $\psi$ for the case of constant appetite (see Theorem 1 in \cite{HHP1}) that  can be easily extended to the general case of random appetites. Such construction uses the Gale-Shapley algorithm (see \cite{GS}) and we will present it here for the sake of completeness considering our case of random appetites: 

\textbf{Gale-Shapley algorithm.} For each positive integer $n$, the stage $n$ consists of the two following parts:

\begin{itemize}
\item[(a)] Each site $x$ non-equidistant from two or more centers \textit{applies} to the closest center to $x$ which has not rejected $x$ at any earlier stage (Note that the set of sites equidistant from at least two centers is a $\mathcal L$-null set. Then, for convenience we can consider it a subset of $\psi^{-1}(\Delta)$.)
\item[(b)] For each center $\xi$, let $A_n(\xi)$ be the set of sites which applied to $\xi$ in the stage $n$ (a), and define the \textit{rejection radius}
\begin{align*}
	r_n(\xi)=\inf \{r:\mathcal L(A_n(\xi)\cap B(\xi,r))\geq\alpha V_\xi\},
\end{align*}
where $B(\xi,r)$ is the open ball with radius $r$ centered at $\xi$, and the infimum of the empty set is taken to be $\infty$. Then $\xi$ shortlists all sites in $A_n(\xi)\cap B(\xi, r_n(\xi))$ and rejects all sites in $A_n(\xi)\setminus B(\xi, r_n(\xi))$. 
\end{itemize}
Consider a site $x$ non-equidistant from two or more centers. Either $x$ is rejected by every center (in order of increasing distance from $x$) or for some center $\xi$ and some stage $n$, $x$ is shortlisted by $\xi$ at every stage after $n$. In the former case, we put $\psi(x)=\infty$ ($x$ is unclaimed); in the latter case, we put $\psi(x)=\xi$.
\\

In Section 2 we will mention some results showing that, under some conditions, this is the only construction with the property of stability. Now, in \cite{Diazldsm} we presented a formal construction of the probability space, based on the construction proposed by Meester and Roy for the Boolean model of continuum percolation (see \cite{MR}, pp. 16--17). We will not present it here, but suffice it to say that it entails the independence between appetites and centers, and the independence among the appetites for different centers. Besides, it also entails the application of some ergodic properties needed below. 

We denote the model by $(\Xi,\alpha V,\psi)$ with joint law $\mathcal P$ and expectation $\mathcal E$. Since $\psi$ is \textit{translation-equivariant} (i.e., $\psi_\Xi(x)=\xi$ implies that $\psi_{T\Xi}(Tx)=T\xi$, for every translation $T\in\mathbb R^d$), we have that $\mathcal P$ is translation invariant. Its \textit{Palm version} $(\Xi^*,\alpha V^*,\psi^*)$ with law $\mathcal P^*$ and expectation operator $\mathcal E^*$ is conditioned to have a center on the origin. Note that since $\alpha V$ is independent of $\Xi$, then $\alpha V^*$ and $\Xi^*$ are also independent and $\alpha V=\alpha V^*$. 

Thorisson showed in \cite{Thor} that $\Xi$ and $\Xi^*$ can be coupled so that almost surely one is a translation of the other; i.e., $\Xi^*=\Xi-\psi^{-1}(0)$ (namely, centering the plane at the Poisson point to which the origin is allocated in the stable marriage gives the Palm version of the Poisson process).  Therefore, $\psi_{\Xi^*}$ is defined $\mathcal P^*$-a.s. Finally, remember that since $\Xi$ is a Poisson process, $\Xi^*$ is a Poisson process with an added point at the origin. (See a similar reasoning in \cite{HHP1} p. 1256 and \cite{HPPS} pp. 8--9). 

\section{Previous results}

Besides the existence of the stable allocation for a discrete set of centers, in \cite{Diazphd} there were some generalizations of previous results in \cite{HHP1}. Here we mention, without proof, some relevant versions for our purpose:

\begin{teo}[\textbf{Almost sure uniqueness}]\label{unique}
	Let $\Xi$ be a Poisson point process with law \textbf P and finite intensity $\lambda$, where each 		center $\xi$ chooses independently its appetite with law $\mathbb P$, then there exists an $\mathcal 	L$-a.e.\ unique stable allocation $\psi$ to $\Xi$, $\mathcal P$-a.s.
\end{teo}

\begin{teo}[\textbf{Phase transitions}]\label{phases}
	Let $\Xi$ be a Poisson point process with law \textbf P and intensity $\lambda\in(0,\infty)$.
	\begin{itemize}
		\item If $\lambda\alpha\mathbb EV<1$ (\textbf{subcritical}), then all the centers are sated but 				there exists an infinite volume of unclaimed sites, $\mathcal P$-a.s..
		\item If $\lambda\alpha\mathbb EV=1$ (\textbf{critical}), then all the centers are sated and $				\mathcal L$-a.e.\ site is claimed, $\mathcal P$-a.s.
		\item If $\lambda\alpha\mathbb EV>1$ (\textbf{supercritical}), then there are unsated centers 				but $	\mathcal L$-a.e.\ site is claimed, $\mathcal P$-a.s.
	\end{itemize}
\end{teo}

\begin{prop}[\textbf{Monotonicity on the centers}]\label{monocenters}
	Let $\Xi_1, \Xi_2$ be two sets of centers chosen according to Poisson processes such that $			\Xi_1\subseteq\Xi_2$. Take $\{V_\xi\}_{\xi\in\Xi_2}$, a set of \textit{iid} r.v.\ distributed like $V$. Fix $	\{v_\xi\}_{\xi\in\Xi_2}$ for a realization of $\{V_\xi\}_{\xi\in\Xi_2}$ as the appetites of the centers $		\Xi_2$, and take $\{v_\xi:\xi\in\Xi_1\cap\Xi_2\}$ for a realization of the appetites of  $\Xi_1$. Let $		\psi_1$ and $\psi_2$ be their respective ($\mathcal L$-a.e.\ unique) stable allocations, then either $	\psi_1(x)=\Delta$ or $|x-\psi_1(x)|\geq|x-\psi_2(x)|$ for $\mathcal L$-a.e.\ $x\in\mathbb R^d$, $\mathcal P$-a.s.
\end{prop}

\begin{prop}[\textbf{Monotonicity on the appetites}]\label{monoappet}
	Let $\Xi$ be a set of centers, let $\{v_\xi\}_{\xi\in\Xi}$ be a realization of $\{V_\xi\}_{\xi\in\Xi}$. Let $		\psi_1$ and $\psi_2$ be two stable allocations with appetites $\alpha_1v_\xi$ and $\alpha_2v_\xi		$, respectively, where $\alpha_1\leq\alpha_2$, then $|x-\psi_1(x)|\geq|x-\psi_2(x)|$, for $\mathcal L	$-a.e.\ $x\in\mathbb R^d$, $\mathcal P$-a.s.
\end{prop}

\textbf{Remark 1}. Theorems \ref{unique} and \ref{phases} have more general statements in \cite{HHP1} (and \cite{Diazphd}). Theorem \ref{unique} was proved for every translation-invariant point process, and Theorem \ref{phases} was proved for every point process ergodic under translations. For our interests, the versions here stated will suffice. 

\textbf{Remark 2}. For the two last propositions, $\mathcal C_1\subseteq\mathcal C_2$, where $\mathcal C_1$ and $\mathcal C_2$ are the closure of the claimed sites for the stable allocations $\psi_1$ and $\psi_2$, respectively.

\textbf{Remark 3}. Proposition \ref{monoappet} admits another way to generalize it (also found in \cite{Diazphd}): Let $V_1$ and $V_2$ be two independent nonnegative random variables such that $V_1$ is stochastically dominated by $V_2$. Then there exists a coupling $(\hat V_1,\hat V_2)$ of $V_1$ and $V_2$ such that $\mathcal C_1\subseteq\mathcal C_2$, for the stable allocations $\psi_1$ and $\psi_2$, respectively, where the appetites of the centers are chosen according to the law of $\hat V_1$ for the allocation $\psi_1$, and with the law of $\hat V_2$ for the allocation $\psi_2$.

\section{Statement of results}

\begin{defi}[\textbf{Percolation}]\label{defperco}
	There is percolation if there exists an unbounded connected subset of $\mathcal C$.
\end{defi}

Let us denote by $A_\alpha$ the event $\{$0 belongs to a connected unbounded subset of $\mathcal C$ in the $d$-dimensional model with appetites given by $\alpha V\}$.

\begin{defi}[\textbf{Critical point}]\label{crit}
	\begin{align*}
		\alpha_p(d):=\sup\{\alpha:\mathcal P[A_\alpha]=0\}.
	\end{align*}
\end{defi}

It is noted immediately that, by Theorem \ref{phases}, there is percolation of $\mathcal C$ whenever $\alpha\geq(\lambda\mathbb EV_\xi)^{-1}$, since $\mathcal L$-a.e.\ site is claimed, $\mathcal P$-a.s. Therefore, one interesting question is whether there exists $\alpha_p(d)>0$ (i.e., whether there exits a positive $\alpha$ such that, almost surely, there is no percolation of $\mathcal C$). The result to be proved in this article is:

\begin{teo}\label{Perco2}
	(i) Let $\{V_\xi\}$ be a collection of nonnegative \textit{iid} r.v.'s distributed like $V$; suppose that 		there exists  $\delta\in\mathbb R^+$ such that $\int_{v\geq0}v^{2+\delta}d\mathbb P(v)<\infty$. Let 		$D$ be the diameter of the connected component of $\mathcal C$ containing the origin. Then, for 		every $s$, there is a small enough $\alpha$ such that $\mathcal ED^s<\infty$. Thus, there is not 			percolation of $\mathcal C$.
	
	(ii) Under the same conditions, for small enough $\alpha$ there is percolation of $\mathbb R^d		\setminus \mathcal C$ (i.e., there is an unbounded connected component of the unclaimed sites).
\end{teo}

Take $\delta_1\in\mathbb R^+:=(0,\infty)$ and consider a random variable $V'=\max\{V,\delta_1\}$; call its law $\mathbb P'$ and its expectation $\mathbb E'$; then $V$ is stochastically dominated by $V'$. Denote by $\mathcal P'$ and by $\mathcal E'$ the probability and the expectation, respectively, of the model $(\Xi, \alpha V',\psi)$ and call $\mathcal C'$ the closure of the set of its claimed sites. 

The r.v.\ $V'$ is chosen because of a proof technicality. However, given Remark 3 and Theorem \ref{Perco2}, there exists a coupling $(V,V')$ of $V$ and $V'$ such that there is not percolation of the claimed sites by the model with appetites given by $\alpha V$ for small enough $\alpha$.

The restriction of the integral $\int_{v\geq0}v^{2+\delta}d\mathbb P(v)$ to be finite is needed to apply Lemma \ref{Deviations}, a result by Nagaev for large deviations of sums of random variables (more on this on Subsection \ref{largedevi}). Now, Lemma \ref{Deviations} was not the only result on this respect proposed by Nagaev in \cite{Nag}; he has some results requiring less in these ``moments'', but we were unable to apply them successfully into this proof. 

It is worth noting that the first result known about the absence of percolation for the stable marriage of Poisson and Lebesgue was due to Freire, Popov and Vachkovskaia (\cite{FPV}). It was obtained for the case in which the appetites are constant, using methods from fractal percolation. In order to prove Theorem \ref{Perco2}, we will rely on a strong result by Gou\'er\'e on subcritical regimes of Boolean models. It was based on previous developments by the same author (see \cite{Gouere1}) in which he improved the phase transition results for the Poisson Boolean model (see \cite{MR}). Before stating it, we define some concepts.

Let $\Pi$ be a point process on $\mathbb R^d$.  The points of $\Pi$ will be centers of balls with i.i.d.\ random radii sampled from $\rho$, a r.v.\ with law $\nu$ on $(0,\infty)$. The \textit{Boolean model} can now be defined: The points of $\Pi$ will be centers of disks such that the marginal distribution of each radius must be the same as $\rho$. Denote this model by $(\Pi,\rho)$.  Note that $(\Pi,\rho)$ can indeed be seen as a point process $\tilde\Pi$ on $\mathbb R^d\times(0,\infty)$. Call the joint probability of the Boolean model $P_\nu$, with associated expectation $E_\nu$.

\begin{teo}\label{Gouere}
 	 Let $S$ be the connected component  of the Boolean model containing the origin. Also, define 		$M:=\sup_{x\in S}\|x\|$ and let $C>0$. Now consider the following four properties:\\
	
	\noindent \textbf{P1.} For all $r>0$ and all $x\in\mathbb{R}^d\setminus B(0, Cr)$, these two 			point process are independent:
		\begin{equation*}
			\tilde\Pi\cap B(0,r)\times(0,r]\text{ and }\tilde\Pi\cap B(x,r)\times(0,r]
		\end{equation*}
		
	\noindent \textbf{P2.} 
		\begin{equation*}
			\sup_{r>0}r^d\nu([r,\infty))\leq H;
		\end{equation*}
	\noindent \textbf{P3.}
			\begin{equation*}
				\int_{[1,\infty)}\beta^d\nu(d\beta)<\infty;
			\end{equation*}
	
	\noindent \textbf{P4.} For some $s$, a positive real,
			\begin{equation*}
				\int_{[1,\infty)}\beta^{d+s}\nu(d\beta)<\infty;\\
			\end{equation*}	
	where $H$ is  positive constant which depends only on the dimension $d$ and the constant C. 		Then, if properties \textbf{P1--P3} are satisfied, then $S$ is $P_\nu$-a.s.\ bounded. If properties 		\textbf{P1--P4} are satisfied, then $E_\nu [M^s]$ is finite.
\end{teo}

\textbf{Remark 4}. Theorem \ref{Gouere} was originally stated for $\nu$ being a locally finite measure. We have kept its original formulation. However, as Gou\'er\'e notes in \cite{Gouere2}, provided that $\nu$ is a finite measure, conditions \textbf{P2} and \textbf{P3} together are equivalent to $\int_{(0,\infty)}\beta^d\nu(d\beta)<\infty$, and conditions \textbf{P2} and \textbf{P4} together are equivalent to $\int_{(0,\infty)}\beta^{d+s}\nu(d\beta)<\infty$.

The proof of Theorem \ref{Gouere} is quite involved and can be found in \cite{Gouere2}. The next two paragraphs present a sketch to understand the first part of the proof (that $S$ is $P_\nu$-a.s.\ bounded). This sketch will allow us to use some of these concepts for the proof of Theorem \ref{Perco2} \textit{(ii)}.

Let $(\Pi,\rho_\gamma^\beta)$ be the Boolean model restricted to the balls with radii in the closed interval $[\gamma,\beta]$, where $\gamma>0$, $\beta>0$. Say that the event $G(x,\gamma,\beta)$ occurs if and only if there is a connected component from $B(x,\beta)$ to $\mathbb R^d\setminus B(x,2\beta)$ using only balls of $(\Pi,\rho_\gamma^\beta)$.

To prove the first part of Theorem \ref{Gouere} it suffices to show that $P_\nu[M>\beta]\rightarrow0$ when $\beta\rightarrow\infty$. To do this, Gou\'er\'e proves that $P_\nu[M>\beta]$ is bounded above by $p_\nu^G(\beta):=P_\nu[G(0,0,\beta)]$ plus an error term that is bounded and that goes to 0 when $\beta$ goes to $\infty$:

\begin{equation*}
	P_\nu[M>\beta]\leq p_\nu^G(\beta)+\tilde D\int_{[\beta,\infty)}r^d\nu(dr),
\end{equation*}
where $\tilde D$ is a positive constant. However, as a result of Proposition 2.1 in \cite{Gouere2},

\begin{equation}\label{pG}
	p_\nu^G(\beta)\rightarrow0  \text{ when } \beta\rightarrow\infty. 
\end{equation}
See also the Remark and the idea of the proof just below the Proposition 2.1 in \cite{Gouere2} to understand better this convergence, but note that Gou\'er\'e calls $\pi(0,\beta)$ what we have called $p_\nu^G(\beta)$.

\section{Proofs}

The proof of Theorem \ref{Perco2} \textit{(i)} presented here is based on the work by Gou\'er\'e in \cite{Gouere2}. The part \textit{(ii)} is proved with an idea developed by Freire, Popov and Vachkovskaia in \cite{FPV}. However, part \textit{(i)} of Theorem \ref{Perco2}  can be proven using fractal percolation, as in \cite{FPV}; nonetheless, the methodology proposed in \cite{Gouere2} allows a more general result (the finiteness of $\mathcal EM^s$ for all positive $s$) and is simpler than that of \cite{FPV}.\\

Let $\Xi$ be a Poisson point process in $\mathbb R^d$ whose intensity is $\lambda>0$. For every $\xi\in\Xi$ we define $R:=R(\xi,\Xi)$ as follows:
\begin{equation*}
	R(\xi, \Xi):=\inf\left\{r>0:\sum_{\zeta\in B[\xi, 2r]}\alpha V_\zeta'\leq\pi_dr^d\right\},
\end{equation*} 
$\inf\emptyset=\infty$. Here $\pi_dr^d$ is the volume of a $d$-dimensional ball with radius $r$; $B[\xi,2r]$ is the closed ball with center $\xi$ and radius $2r$. Note that there is at least one center in the closed ball $B[\xi,2R]$ (namely $\xi$) and its appetite would be at least $\alpha\delta_1$ (recall that $\delta_1$ is the lower bound for $V'$); therefore, the minimum value of $R$ is $(\alpha\delta_1\pi_d^{-1})^{1/d}$. We will denote this value by $b$.  We can now construct a Boolean model with its centers being $\Xi$ and the radii of every $\xi\in\Xi$ given by $R(\xi,\Xi)$. In agreement with the former notation, we call this model $(\Xi,R)$. Note that $\mathcal P'$ is also the law of $R$, by construction.

The idea behind this concept was first used in \cite{HHP2} to prove some tail bounds; after, it was also used in \cite{FPV} to prove the absence of percolation when the appetites are constant.

\subsection{Domination}

\begin{lema}\label{contem}
Let $\Sigma$ be the union of the balls of the Boolean model previously constructed. If $\alpha$ is small enough, then $\mathcal C'$ is contained in $\Sigma$, $\mathcal P'$-a.s.
\end{lema}

\begin{proof}
	Let $\xi\in\Xi$. Using the strong law of large numbers, it is easy to see that $R(\xi,\Xi)$ is finite with  		probability one, provided 
	
	\begin{equation}\label{alphabounded}
		\alpha< (\lambda2^d\mathbb E'V')^{-1}. 
	\end{equation}
	
	To prove the Lemma we only need to check that $\psi^{-1}(\xi)\subset B[\xi,R]$, where $R=R(\xi,\Xi)$. 	Thus, under (\ref{alphabounded}), we have
	\begin{equation*}
		\sum_{\zeta\in B[\xi, 2R]}\alpha V_\zeta'=\pi_dR^d,	
	\end{equation*}
	for small $\alpha$. Take $\varepsilon>0$ such that there are no points of $\Xi$ in $ B[\xi,2R			+2\varepsilon]\setminus B[\xi,2R]$ (Note that this $\varepsilon$ exists with probability 1). Then
	\begin{equation*}
		\sum_{\zeta\in B[\xi, 2R+2\varepsilon]}\alpha V_\zeta'<\pi_d(R+\varepsilon)^d.
	\end{equation*}
	Therefore, by the definition of territories,
	\begin{equation*}
		\mathcal L\left[\psi^{-1}\left(\Xi\cap B[\xi, 2R+2\varepsilon]\right)\right]<\pi_d(R+\varepsilon)^d.
	\end{equation*}
	Then, there exists $x\in B[\xi,R+\varepsilon]$ such that $\psi(x)\in\Xi\cup\{\infty\}\setminus B[\xi,2R			+2\varepsilon]$. If $\psi(x)\in\Xi$,
	\begin{equation*}
		|x-\psi(x)|>R+\varepsilon\text{ and } |x-\xi|\leq R+\varepsilon.
	\end{equation*}
	Thus, $x$ desires $\xi$. If $\psi(x)=\infty$, $x$, by definition, desires $\xi$. Since $\psi$ is stable, $\xi		$ does not covet $x$. As a consequence, $\psi^{-1}(\xi)\subset B[\xi,|x-\xi|]$ and then $				\psi^{-1}(\xi)\subset B[\xi,R+\varepsilon]$. Taking $\varepsilon>0$ arbitrarily small, we 				obtain that $\psi^{-1}(\xi)\subset B[\xi,R]$.
\end{proof}

\subsection{Large deviations}\label{largedevi}

We need a classical result of Nagaev on large deviations of sums of random variables (Corollary 1.8 of \cite{Nag}) and also the Chernoff bounds for Poisson random variables. Those results will be used to bound the the probability of the event $\{R(\xi,\Xi)>r/2\}$, where $r$ is a constant.

\begin{lema}[\textbf{Large deviations of sums}]\label{Deviations}
	Let $S_n=\sum_{i=1}^nZ_i$, where $\{Z_i\}$ is a sequence of \textit{iid} r.v.\ with law $P$, 0-mean and variance $\sigma^2$. Define $A_{2+\delta}^+:=\int_{u\geq0}u^{2+\delta}dP(u)<\infty$ ($\delta$ as defined on Theorem \ref{Perco2}). Then
	\begin{equation*}
		\text P\left[S_n>x\right]\leq\left(\frac{4+\delta}{2+\delta}\right)^{2+\delta}nA_{2+	\delta}^+x^{-(2+		\delta)}+\exp\left\{\frac{-2e^{-(2+\delta)}x^2}{(4+\delta)^2n\sigma^2}\right\}
	\end{equation*}	
\end{lema}

\textbf{Remark 5}. In \cite{Nag}, Lemma \ref{Deviations} was more general: the r.v.'s did not need to be identically distributed and the integral $\int_{u\geq0}u^tdP(u)$ existed for some $t\geq2$.

\textbf{Remark 6}. Although there are other Nagaev results for moments less than 2, those results are not useful for our purposes, as will be explained below. Indeed, we need something more than the finiteness of the second moments ---we need the finiteness of $A_{2+\delta}^+$.

\begin{lema}[\textbf{Poisson large deviations}]\label{ChePo}
	Let $N$ be a r.v.\ with Poisson distribution and mean $0<\tilde\lambda<\infty$ (denoted $N\sim			\mathcal P_{\tilde\lambda})$, then
	\begin{align*}
		\mathcal P_{\tilde\lambda}[N\geq a]&\leq e^{-\tilde\lambda g(\tilde\lambda/a)}\text{\ \ \ for } a>			\tilde\lambda
	\end{align*}
	where $g(x)=(x-1-\log x)/x$ and $g(x)\rightarrow\infty$ when $x\rightarrow0$.
\end{lema}

\begin{lema}\label{lambda}
	Given the same set of conditions of Theorem \ref{Perco2}, for $r>0$, then $\mathcal P'[R>r]\leq c			r^{-d(1+\delta)}$, where $c\in\mathbb R^+$.
\end{lema}

\begin{proof}
	Let $r>0$ and call $\mathcal P'^*$ the Palm version of $\mathcal P'$. By definition of $R$ and $\Xi$ we have:
	\begin{align*}
		\mathcal P'[R>r]&=\mathcal E'\left(\sum_{\xi\in\Xi\cap[0,1]^d}\textbf 1_{\{R(\xi,\Xi)>r\}}\right)\\
		&\leq\lambda\mathcal P'^*\left[R(\xi,\Xi)>r/2\right]\\
		&\leq\lambda\mathcal P'^*\left[\sum_{\zeta\in B[\xi,2r]}V_\zeta'>\pi_dr^d2^{-d}\alpha^{-1}\right],
	\end{align*}
	where the first inequality is obtained using the definition of expectation in a compound (Poisson) 		point process and because of the fact that $R$ is always positive. Let $\lceil y\rceil$ be the closest 		integer higher than $y$ (or $y$, if $y\in\mathbb Z$). Fix $m\in(2,r)$, and let $K=\left\lceil\lambda		\pi_dm^d(r^d+r^{d/2+\varepsilon})\right\rceil$ and $N\sim\mathcal P(\lambda\pi_dr^dm^d)$, this way 	$N$ dominates $\#\{\Xi\cap B[\xi,2r]\}$. Call $\mu$ the expectation of $V'$ and let $\hat V_i=V_i'-\mu$. 	Finally take $\alpha^{-1}>\tilde m^d\lambda2^d\mu$, where $\tilde m>m$ (note that it satisfies 		(\ref{alphabounded})). With these considerations we have:
	\begin{align}\label{inequ}
		\mathcal P'[R>r]&\leq\lambda\mathcal P'^*\left[\sum_{k=1}^NV_k'>\pi_dr^d2^{-d}\alpha^{-1}\right]			\nonumber\\
		&\leq\lambda\mathcal P'^*\left[\sum_{k=1}^NV_k'>\pi_d\tilde m^d\lambda\mu r^d\right]					\nonumber\\
		&\leq\lambda\left\{\mathcal P'^*\left[\sum_{k=1}^NV_k'>\pi_d\tilde m^d\lambda\mu r^d|N\leq K				\right]\mathcal P'^*[N\leq K]+\mathcal P'^*[N>K]\right\}\nonumber\\
		&\leq\lambda\left\{\mathcal P'^*\left[\sum_{k=1}^K\hat V_k>x\right]+\mathcal P'^*[N>K]\right\},
	\end{align}
	where $x:=\pi_d\tilde m^d\lambda\mu r^d-K\mu=\pi_d\lambda\mu r^d[\tilde m^d-m^d(1+r^{-d/2+		\varepsilon})]$ (note that $x>0$ for large enough $r$, and that $x$ and $K$ have order $r^d$). Thus, 	using Lemma \ref{Deviations} for the first term in the braces in \ref{inequ}, we obtain for large enough 	$r$ that
	\begin{align}\label{PV}
		\mathcal P'^*\left[\sum_{k=1}^K\hat V_k>x\right]\leq c_3r^{-d(1+\delta)}+e^{-c_2r^d},
	\end{align}
	for suitable positive constants $c_3$ and $c_2$. Now we use Lemma \ref{ChePo} for the second term 	in (\ref{inequ}): considering that $g(x)$ in this Lemma has order $(x-1)^2/2$ when $x\rightarrow1$, 	we get that $g(r^d/(r^d+r^{d/2+\varepsilon}))\sim r^{-d+2\varepsilon}$ when $r\rightarrow\infty$. Thus 	$\lambda\pi_dr^dm^dg(r^d/(r^d+r^{d/2+\varepsilon}))\sim r^{2\varepsilon}$ when $r\rightarrow\infty$ 	and we obtain
	\begin{align}\label{PN}	
		P[N>K]\leq c_1r^{2\varepsilon},
	\end{align}
	for some positive constant $c_1$. Finally, using (\ref{inequ}), (\ref{PV}) and (\ref{PN}), we obtain
	\begin{align*}
		\mathcal P'[R>r]&\leq \lambda\left\{c_3r^{-d(1+\delta)}+e^{-c_2r^d}+e^{-c_1r^d}\right\}\\
		&\leq c\lambda r^{-d(1+\delta)},
	\end{align*}
	for a positive $c$. Note that the last expression goes to 0 when $\lambda\rightarrow0$ (we will use 	this fact to prove Theorem \ref{Perco2}).
\end{proof}

\subsection{Proof of Theorem \ref{Perco2}}

Given Lemma \ref{contem}, we only need to check that our Boolean model satisfies the conditions of Theorem \ref{Gouere}. Now we proceed in the same way as in \cite{Gouere2}:

\textbf{P1:} Let $r>0$. For all $\xi\in\Xi$ define	
	\begin{equation*}
		\tilde R(\xi,\Xi)=\inf\left\{s\in[0,r]:\sum_{\zeta\in B[\xi,2s]}\alpha V_\zeta'\leq\pi_d(\xi)s^d\right\}
	\end{equation*}	
And $\tilde R(\xi,\Xi)=r$ if there is no such $s$. Note that, for all $\xi\in\Xi$, if $R(\xi,\Xi)<r$ or $\tilde R(\xi,\Xi)<r$, then $\tilde R(\xi,\Xi)=R(\xi,\Xi)$. So,
	\begin{equation*}
		\{\xi\in\Xi:R(\xi,\Xi)<r\}=\{\xi\in\Xi:\tilde R(\xi,\Xi)<r\}.
	\end{equation*}
Thus, for every $r>0$, we have that $\{\xi\in\Xi\cap B(x,r):R(\xi,\Xi)<r\}$ only depends on $\Xi\cap B(x,3r)$. By the independence of the Poisson process  we obtain that if $x\in\mathbb R^d\setminus B(0,6r)$, then the point processes $\{\xi\in\Xi\cap B(0,r):R(\xi,\Xi)<r\}$ and $\{\xi\in\Xi\cap B(x,r):R(\xi,\Xi)<r\}$ are independent. Therefore, $C=7$ satisfies the property.

\textbf{P2:} Using in the first equality below that $V'$ is lower bounded by $\delta_1$ (therefore $R$ is lower bounded by $b$); and using after Lemma \ref{lambda}, we obtain
	\begin{align*}
		\sup_{r>0}r^d\mathcal P'[R>r]&=\sup_{r\geq b}r^d\mathcal P'[R>r]\\
		&\leq\sup_{r\geq b}r^dc\lambda r^{-d(1+\delta)}\\
		&=c\lambda\sup_{r\geq b}r^{-d\delta}\\
		&=c\lambda b^{-d\delta},
	\end{align*}	
thus, for small enough $\lambda$, we obtain that the $\sup_{r>0}r^d\mathcal P'[R>r]$ is bounded for a given $H(d,C)$, as required. Note that $r^{-d\delta}$ is decreasing, this is why we needed the ``moment'' $A_{2+\delta}^+$. Also note that, since $r^{-d\delta}$ is decreasing, then the condition $V'>\delta_1$ makes a supremum of the last expression. 
	
\textbf{P3} and \textbf{P4}. By Lemma \ref{lambda} and Remark 4 we obtain: $\int_{(0,\infty)}r^{d+s}d\mathcal P'<\infty$ for every $s\in(0,d\delta]$.

Therefore, given small enough $\lambda$, we can use Theorem \ref{Gouere} to obtain that $\mathcal EM^s<\infty$, and by Lemma \ref{contem}, we get that $\mathcal ED^s<\infty$. Thus, after rescaling,  Theorem \ref{Perco2} \textit{(i)} is proved for fixed $\lambda$ and small enough $\alpha$.\\

To prove Theorem \ref{Perco2} \textit{(ii)} consider the following sequence of events:
\begin{align*}
	&\{\text{every connected subset in } \mathbb R^d\setminus\mathcal C' \text{ is bounded}\}\\
	\subset&\{\text{for any bounded } K\subset\mathbb R^2 \text{ there is a contour around } K \text{ in }			\mathcal C'\cap\mathbb R^2\}\\
	\subset&\{\text{for any bounded } K\subset\mathbb R^2 \text{ there is a contour around } K \text{ in }			(\Xi,R)\cap\mathbb R^2\}\\
	\subset&\{\text{for infinite values of $\beta$, } G(0,0,\beta)\text{ occurs}\},\\
	\subset&\{\text{for infinite values of $\beta$, } G(0,0,\lceil\beta\rceil)\text{ occurs}\}\\
	\subset&\{\text{for infinite values of $n\in\mathbb N$, } G(0,0,n)\text{ occurs}\}
\end{align*}
where $(\Xi,R)$ is the Boolean model with centers given by $\Xi$ and radii given by $R$, and the former to last event contains the former one because of Proposition \ref{monoappet} and Remark 3. By the Borel-Cantelli lemma, the last event has probability 0 when $\alpha$ is small enough. Thus, Theorem \ref{Perco2} \textit{(ii)} is proved.

\section{Acknowledgments}

The author is grateful to the referees for their valuable suggestions and acknowledges the support of CAPES during his research. He wants to thank Carolina Hern\'andez and Juan S\'aenz for their independent proof-reading and style suggestions, he also thanks Ang\'elica Mar\'ia Vega for her timely and valuable comments.

\end{document}